# ESTIMATION IN DIRICHLET RANDOM EFFECTS MODELS

BY MINJUNG KYUNG[1], JEFF GILL[1] AND GEORGE CASELLA[2]

*University of Florida, Washington University and University of Florida*

We develop a new Gibbs sampler for a linear mixed model with a Dirichlet process random effect term, which is easily extended to a generalized linear mixed model with a probit link function. Our Gibbs sampler exploits the properties of the multinomial and Dirichlet distributions, and is shown to be an improvement, in terms of operator norm and efficiency, over other commonly used MCMC algorithms. We also investigate methods for the estimation of the precision parameter of the Dirichlet process, finding that maximum likelihood may not be desirable, but a posterior mode is a reasonable approach. Examples are given to show how these models perform on real data. Our results complement both the theoretical basis of the Dirichlet process nonparametric prior and the computational work that has been done to date.

**1. Introduction.** Linear and generalized linear mixed models have become important statistical tools bringing more flexibility through the addition of random effects to the more traditional linear models. A general mixed effects model can be written as

$$(Y_1, \ldots, Y_n) \sim f(y_1, \ldots, y_n | \theta, \psi_1, \ldots, \psi_n) = \prod_i f(y_i | \theta, \psi_i),$$

$$\psi_i \sim G, \qquad i = 1, \ldots, n,$$

where $f$ and $G$ are often taken to be normal. The addition of a link function turns this into a generalized linear mixed model; we will discuss that below. The restriction of $G$ to a normal distribution is sometimes thought to be limiting. For example, researchers in the social sciences are uncomfortable with this assumption, wanting a more flexible, possibly nonparametric

Received February 2009; revised June 2009.

[1]Supported by NSF Grants DMS-06-31632 and SES-0631588.

[2]Supported by NSF Grants DMS-04-05543, DMS-06-31632 and SES-0631588.

structure here [Gill and Casella (2009)]. Another troublesome fact, as noted by Burr and Doss (2005), is that random effects, unlike error terms, cannot be checked (there are no residuals). Thus we are totally dependent on this uncheckable model assumption.

One way of relaxing this assumption is with a richer, nonparametric model for $\psi$ with a popular alternative being the Dirichlet process

$$\psi_i \sim \mathcal{DP}(m, \phi_0), \qquad i = 1, \ldots, n,$$

where $\mathcal{DP}$ is the Dirichlet process with base measure $\phi_0$ and precision parameter $m$. By moving to this model we not only relax the normal assumption, but also provide a richer model for the random effects. Such a model has potential for capturing more types of variability in those effects with the possible end result of more precise estimates of the fixed effects. Here we will look at ways to fit such models, focusing on methods of estimating the precision parameter $m$, and evaluating and improving on Gibbs samplers for the models.

1.1. *Background.* Dirichlet process mixture models were introduced by Ferguson (1973) who defined the process and investigated basic properties. Blackwell and MacQueen (1973) showed that the marginal distribution of the Dirichlet process is equal to the distribution of the $n$th step of a Pólya urn process. In particular, they proved that for $\psi_1, \ldots, \psi_n$ i.i.d. from $G \sim \mathcal{DP}$, the joint distribution of $\boldsymbol{\psi}$ is a product of successive conditional distributions of the form

$$(1) \quad \psi_i | \psi_1, \ldots, \psi_{i-1}, m \sim \frac{m}{i-1+m} \phi_0(\psi_i) + \frac{1}{i-1+m} \sum_{l=1}^{i-1} \delta(\psi_l = \psi_i),$$

where $\delta$ denotes the Dirac delta function.

Other work that characterizes the properties of the Dirichlet process includes Korwar and Hollander (1973) and Sethuraman (1994). Work that has particular importance for our development is that of Lo (1984), who derives the analytic form of a Bayesian density estimation that is generated by convoluting a known density kernel with a Dirichlet process, and Liu (1996), who derives an identity for the profile likelihood estimator of $m$.

The implementation of the Dirichlet process mixture model has been made feasible by modern methods of Bayesian computation and efficient algorithms. The work of Escobar and West (1995) and MacEachern and Muller (1998) developed estimation techniques and sampling algorithms, and Neal (2000) provided an extended and more efficient Gibbs sampler.

Note that the representation (1) induces clusters in the random effects since with positive probability the value of $\psi_i$ is equal to one of the previous values. McCullagh and Yang (2006) showed that the marginal distribution of these Dirichlet clusters can be derived using cycles of integers and exchangeability based on Pitman (1996), and an exchangeable cluster process can be generated by a standard Dirichlet allocation scheme. Quintana and Iglesias (2003) show that the joint marginal distribution of the Dirichlet observations can be expressed as a product partition model. Such models were introduced by Hartigan (1990) and Barry and Hartigan (1992), and are based on modeling random partitions of the sample space.

It is interesting to note that with $n$ observations from a Dirichlet process with precision parameter $m$, the marginal distribution of a partition $\{n_1, n_2, \ldots, n_k\}$, where $\sum_j n_j = n, n_j \geq 1$, is given by

$$(2) \qquad \pi(n_1, n_2, \ldots, n_k) = \frac{\Gamma(m)}{\Gamma(m+n)} m^k \prod_{j=1}^{k} \Gamma(n_j),$$

which is a normalized probability distribution on the set of all partitions of $n$ observations. This is the same distribution derived by McCullagh and Yang (2006) using cycles of integers (which can be related to the partitions) and has been used by other authors [Crowley (1997), Booth, Casella and Hobert (2008)] as a prior distribution on clusters for a Bayesian clustering algorithm.

Most of theoretical and computational work for the Dirichlet process mixture models focus on the efficient estimation of $\phi_0$. For the Dirichlet prior, we also need to consider the precision parameter $m$ because it strongly influences the number of distinct components, which is the distribution of the underlying random effects. The estimation of $m$ has had difficulties in implementation due to computational intractability. The number of distinct components is not known in practice and it can be changed if new data are observed. Doss (2008) notes that the precision parameter is typically the most difficult to estimate or defend as a fixed value.

The approach of Liu (1996), for the estimation of $m$, is to use sequential imputation to estimate $m$ using maximum likelihood and treating the subclusters as missing data. Naskar and Das (2004, 2006) estimated $m$ using Monte Carlo EM but did not investigate the properties of the solution. Our approach is based on using marginal or profile likelihood for $m$, and using the Gibbs sampler to estimate the model parameters. This is a variation of Casella (2001) who showed that, using an empirical Bayes approach, the hyperparameters can be estimated in a computationally feasible way.

1.2. *Summary.* In this paper, for a mixed Dirichlet random effects model, we develop algorithms for estimation of the precision parameter and MCMC algorithms for fitting the models. We focus on linear models but also show how to extend our results to a generalized Dirichlet process mixed model with a probit link function. In Section 2 we develop our model and its probit extension using a new parameterization of the Dirichlet mixed model. Section 3 shows how to estimate the precision parameter, $m$. There we see that there is, in fact, not much information in the data about $m$ (as it only depends on the size of the subclusters). We find that the MLE of Liu (1996) can be unstable and show how to obtain a more stable posterior mode estimate. Section 4 derives a Gibbs sampler for the model parameters and the subclusters of the Dirichlet process, and we use our new parameterization of the hierarchical model to derive a new Gibbs sampler that more fully exploits the structure of the model and mixes very well. We can adapt the results of Hobert and Marchev (2008) to establish that our sampler is an improvement, in terms of operator norm and efficiency, over other commonly used algorithms. Section 6 given details for the estimation of the precision parameter, and Section 7 contains illustrative applications. Section 8 summarizes these contributions and adds some perspective, and there are a number of technical appendices.

**2. Models and likelihoods.** A general random effects Dirichlet model can be written as

$$(3) \qquad (Y_1, \ldots, Y_n) \sim f(y_1, \ldots, y_n | \theta, \psi_1, \ldots, \psi_n) = \prod_i f(y_i | \theta, \psi_i),$$

$$\psi_i \sim \mathcal{DP}(m, \phi_0), \qquad i = 1, \ldots, n,$$

where $\mathcal{DP}$ is the Dirichlet process with base measure $\phi_0$ and precision parameter $m$. The vector $\theta$ contains all of the model parameters which we will address a bit later.

Applying the Blackwell–MacQueen formula (1), we can calculate the likelihood function, which by definition is integrated over the random effects, as

$$L(\theta | \mathbf{y}) = \int f(y_1, \ldots, y_n | \theta, \psi_1, \ldots, \psi_n) \pi(\psi_1, \ldots, \psi_n) \, d\psi_1 \cdots d\psi_n,$$

where

$$(4) \qquad \pi(\psi_1, \ldots, \psi_n) = \prod_{i=1}^n \frac{m \phi_0(\psi_i) + \sum_{j=1}^{i-1} I(\psi_j = \psi_i)}{m + i - 1}.$$

From Lo's (1984) Lemma 2 and following the development in Liu's (1996) Theorem 1, we can evaluate this integral to get

$$L(\theta | \mathbf{y}) = \frac{\Gamma(m)}{\Gamma(m+n)} \sum_{k=1}^n m^k \sum_{C:|C|=k} \prod_{j=1}^k \Gamma(n_j) \int f(\mathbf{y}_{(j)} | \theta, \psi_j) \phi_0(\psi_j) \, d\psi_j,$$

where $C$ defines the subclusters, $\mathbf{y}_{(j)}$ is the vector of $y_i$s that are in subcluster $j$ and $\psi_j$ is the common parameter for that subcluster. There are $\mathcal{S}_{n,k}$ different partitions $C$, the Stirling number of the second kind.

A partition $C$ clusters of the sample of size $n$ into $k$ groups, $k = 1, \ldots, n$, and we call these "subclusters" since the grouping is done nonparametrically rather than on substantive criteria. That is, it is likely that any real underlying clusters would be broken up into multiple subclusters by the nonparametric fit since there is little penalty for over-separation. Recall that the regular GLMM assumes that the random effect $\psi_i$s are independent and identically distributed with the normal distribution $N(0, \sigma_\psi^2)$. However, the subclustering assigns different normal parameters across groups and the same parameters within groups; cases are i.i.d. only if they are assigned to the same subcluster.

2.1. *A matrix representation of subclusters.* Each partition $C$ can be associated with an $n \times k$ matrix $A$ defined by $A' = (a_1', a_2', \ldots, a_n')$ where $a_i$ is a $1 \times k$ vector of all zeros except for a 1 in one position that indicates which group the observation is from. Note that the column sums of $A$ are $(n_1, n_2, \ldots, n_k)$, the number of observations in the groups, and there are, of course, $\mathcal{S}_{n,k}$ such matrices. Specifically, if the partition $C$ has groups $\{S_1, \ldots, S_k\}$, then if $i \in S_j$, $\psi_i = \eta_j$ and the random effect can be rewritten as $\boldsymbol{\psi} = A\boldsymbol{\eta}$. For example, if $S_1 = \{3, 4, 6\}$, $S_2 = \{1, 2\}$, $S_3 = \{5\}$,

$$(5) \qquad \begin{pmatrix} \psi_1 \\ \psi_2 \\ \vdots \\ \psi_6 \end{pmatrix} = \begin{pmatrix} 0 & 1 & 0 \\ 0 & 1 & 0 \\ 1 & 0 & 0 \\ 1 & 0 & 0 \\ 0 & 0 & 1 \\ 1 & 0 & 0 \end{pmatrix} \begin{pmatrix} \eta_1 \\ \eta_2 \\ \eta_3 \end{pmatrix}.$$

We then have

$$\sum_{C:|C|=k} \prod_{j=1}^{k} \Gamma(n_j) \int f(\mathbf{y}_{(j)} | \theta, \psi_j) \phi_0(\psi_j) \, d\psi_j$$

$$= \sum_{A \in \mathcal{A}_k} \prod_{j=1}^{k} \Gamma(n_j) \int f(\mathbf{y} | \theta, A\boldsymbol{\eta}) \phi_0(\boldsymbol{\eta}) \, d\boldsymbol{\eta},$$

where $\mathcal{A}_k$ is the set of all matrices $A$ and $\eta_j \sim \phi_0$, independent. If we define

$$(6) \qquad f(\mathbf{y} | \theta, A) = \int f(\mathbf{y} | \theta, A\boldsymbol{\eta}) \phi_0(\boldsymbol{\eta}) \, d\boldsymbol{\eta},$$

the likelihood function is

$$(7) \qquad L(\theta | \mathbf{y}) = \frac{\Gamma(m)}{\Gamma(m+n)} \sum_{k=1}^{n} m^k \sum_{A \in \mathcal{A}_k} \prod_{j=1}^{k} \Gamma(n_j) f(\mathbf{y} | \theta, A).$$

Note that if the integral in (6) can be done analytically, as will happen in the normal case discussed next, we have effectively eliminated the random effects from the likelihood, replacing them with the $A$ matrices, which serve to group the observations.

2.2. *Marginalizing the random effects.* Start with a normal linear model,

$$\mathbf{Y}|\psi \sim N(X\beta + \boldsymbol{\psi}, \sigma^2 I),$$

where $\boldsymbol{\psi} = (\psi_1, \ldots, \psi_n)'$, $\psi_i \sim \mathcal{DP}(m, N(0, \tau^2))$, $i = 1, \ldots, n$, independent. When we introduce the $A$ matrices we get the models

$$\mathbf{Y}|\boldsymbol{\eta}, A \sim N(X\beta + A\boldsymbol{\eta}, \sigma^2 I), \qquad \boldsymbol{\eta} \sim N_k(0, \tau^2 I),$$

and marginally, $\mathbf{Y}$ is multivariate normal with

$$\mathrm{E}\mathbf{Y} = X\beta, \qquad \mathrm{Var}\,\mathbf{Y} = \sigma^2 I + \tau^2 AA'.$$

Moreover, in this model we can analytically marginalize out some of the model parameters. See Section 4.

2.3. *Consistency.* We briefly discuss posterior consistency of these models, noting that they have been extensively examined by Ghosal and co-authors [see, in particular, Ghosal, Ghosh and Ramamoorthi (1999) and Ghosal (2009)]. Ghosal (2009) defines a *Mixture of Dirichlet Process* (MDP), where the Dirichlet process is the error distribution, with possibly additional hyperparameters for the base measure. Alternatively, there is the *Dirichlet Process Mixture* (DPM)[3] where, for a parametric density, we assume there is a latent variable from a unknown and unrestricted distribution, modeled with a Dirichlet process prior. The famous inconsistency result of Diaconis and Freedman (1986), and the results of Doss (1985a, 1985b), are for the MDP model while our Dirichlet random effects model is an example of DPM with a normal density for the observations.

Ghosal (2009) showed that for the MDP model, the conditions of the general consistency theory of Schwartz (1965) are not satisfied due to the discreteness of the resulting Dirichlet likelihood; thus the posterior from the MDP model can be inconsistent. However, the DPM likelihood is smooth enough to satisfy the conditions of Schwartz's theorem, and thus posterior consistency holds. Therefore, it follows that the posterior of our proposed Dirichlet random effects model is consistent.

We conducted a small simulation study to examine the behavior of the posterior estimates of the coefficients in a linear DPM model. In this simulation, we fix the number of subclusters, $k$, to be 20% of sample size, and the concentration parameter $m$ by using (18). Other parameters are fixed at $\sigma^2 = 1$, $\tau^2 = 4$ and $\boldsymbol{\beta} = (\beta_0, \beta_1)' = (1, 2)'$. From Table 1 we see that the standard deviations get smaller as sample size bigger, and the estimates are getting close to the true value as $n$ gets bigger, illustrating the posterior consistency for the linear Dirichlet random effects models.

---

[3]It appears that the terms *Mixture of Dirichlet Process* and *Dirichlet Process Mixture* and not used unambiguously in the literature.

2.4. *A probit mixed Dirichlet random effects model.* A generalized linear *mixed* model (GLMM) can be specified to accommodate outcome variables conditional on mixtures of possibly correlated random and fixed effects [Breslow and Clayton (1993)]. For example, assume that there is a Bernoulli selection process where we observe $Y_i$ according to

$$(8) \qquad Y_i \sim \text{Bernoulli}(p_i), \qquad i = 1, \ldots, n,$$

where $y_i$ is 1 or 0; thus $p_i = \text{E}(Y_i)$ is the probability of a success for the $i$th observation. Moreover, using a link function $g(\cdot)$, we can express the transformed mean of $Y_i$ as a linear function,

$$(9) \qquad g(p_i) = \mathbf{X}_i \boldsymbol{\beta} + \psi_i,$$

where $\mathbf{X}_i$ are covariates associated with the $i$th observation, $\boldsymbol{\beta}$ is the coefficient vector, and $\psi_i$ is a random effect accounting for subject-specific deviation from the underlying model. The $\psi_i$s are usually assumed to be distributed as $N(0, \sigma^2_\psi)$.

Variations of GLMMs were used by Dorazio et al. (2008) to model spatial heterogeneity in animal abundance, and Gill and Casella (2009) modeled political science data by using a GLMM with an ordered probit link. For Bayesian inference, Albert and Chib (1993) used a Gibbs sampler by introducing a latent variable into the model and noted that the probit model on the binary response is connected with the normal linear model on a continuous latent data response. Mukhopadhyay and Gelfand (1997) used a fully Bayesian approach to fit generalized linear models with Dirichlet random effects, and Ghosh et al. (1998) proposed hierarchical Bayes generalized linear models for longitudinal GLMMs in small area estimation, and provided a general theorem that ensures the propriety of posteriors under diffuse priors.

TABLE 1
*Estimation of the coefficient parameters: Posterior means and standard deviations*

| Condition | $m$ | $\beta_0$ | $\beta_1$ |
|---|---|---|---|
| $n = 100$ and $k = 20$ | 7.246 | 0.618 (0.114) | 1.945 (0.110) |
| $n = 500$ and $k = 100$ | 37.318 | 0.927 (0.046) | 1.986 (0.045) |
| $n = 1000$ and $k = 200$ | 74.915 | 1.051 (0.032) | 2.040 (0.033) |

Although we can proceed to develop estimators with a general link function (9), choosing $g$ to be the probit link greatly simplifies things, which is what we will do. If we introduce the latent variable

$$U_i \sim N(X_i\beta + \psi_i, \sigma^2), \tag{10}$$

then defining $Y_i = I(U_i > 0)$ results in the probit model. Thus, if we consider the hierarchy to start with the $U_i$, we are in exactly the model of Section 2.2. To fit the probit model we then add an additional step to the Gibbs sampler to generate the $U_i$ conditional on the remaining parameters. This results in a truncated normal random variable generation, with details in Appendix A.2.

**3. Estimating the precision parameter.** In this section we look at the performance of the maximum likelihood estimates of $m$ using both profile likelihood and marginal likelihood. We find that these estimates can be unstable, and hence suggest an alternative based on a posterior mode, which we ultimately implement with importance sampling in Section 6.

3.1. *Maximum likelihood estimates.* From (7) define

$$\mathcal{L}_k(\theta|\mathbf{y}) = \sum_{A \in \mathcal{A}_k} \prod_{j=1}^{k} \Gamma(n_j) f(\mathbf{y}|\theta, A). \tag{11}$$

Then, as in Liu (1996) or Doss (1994), we can obtain a *profile likelihood estimate* of $m$. That is, for each fixed value of $\theta$ we can differentiate the log of (7) and set it equal to 0 to obtain the profile MLE as the solution to

$$\frac{\sum_{k=1}^{n} k m^{k-1} \mathcal{L}_k(\theta|\mathbf{y})}{\sum_{k=1}^{n} m^k \mathcal{L}_k(\theta|\mathbf{y})} = \sum_{i=1}^{n} \frac{1}{m+i-1} \tag{12}$$

which defines $\hat{m}(\boldsymbol{\theta})$, the profile MLE.

This development can be easily modified to obtain the marginal likelihood estimate of $m$, $\hat{m}$. Under the usual regularity conditions, when $\hat{m}$ and $\hat{m}(\hat{\theta})$ are both consistent estimates, we would expect that the marginal MLE would be a more stable estimate of the precision parameter.

Starting from (11) we can also integrate over $\theta$ to obtain the marginal likelihood,

$$\mathcal{L}_k(\mathbf{y}) = \int \mathcal{L}_k(\theta|\mathbf{y}) \, d\theta. \tag{13}$$

The same development as above will lead to the expression (12) with $\mathcal{L}_k(\theta|\mathbf{y})$ replaced by $\mathcal{L}_k(\mathbf{y})$.

We also note that expressions for the approximate variances of either $\hat{m}$ or $\hat{m}(\hat{\theta})$ can be easily attained from the second derivative of the appropriate log likelihood evaluated at the estimate.

3.2. *Solving the likelihood equations.* Note that, for either profile or marginal likelihood, we can write the log likelihood function as

$$(14) \qquad \ell(m) = \log\left(\sum_{k=1}^{n} m^k c_k\right) - \sum_{i=1}^{n} \log(i - 1 + m),$$

where $c_k = \mathcal{L}_k(\theta|\mathbf{y})$ for profile likelihood and $c_k = \mathcal{L}_k(\mathbf{y})$ for marginal likelihood. The derivative of the log-likelihood is given by

$$(15) \qquad \frac{\partial}{\partial m}\ell(m) = \frac{\sum_{k=1}^{n} km^{k-1}c_k}{\sum_{k=1}^{n} m^k c_k} - \sum_{i=1}^{n} \frac{1}{i - 1 + m}.$$

It is straightforward to show that

$$(16) \qquad \lim_{m \to 0} \frac{\partial}{\partial m}\ell(m) = \frac{2c_2}{c_1} - \sum_{i=2}^{n} \frac{1}{i - 1},$$

$$\lim_{m \to \infty} \mathrm{sgn}\left(\frac{\partial}{\partial m}\ell(m)\right) = \mathrm{sgn}\left(\sum_{i=1}^{n}(i-1) - \frac{c_{n-1}}{c_n}\right),$$

where $\mathrm{sgn}(\cdot)$ is the sign of the function. Note that $\lim_{m \to \infty} \frac{\partial}{\partial m}\ell(m) = 0$, but the direction of approach is important. The signs of the derivatives at the extremes are only functions of the ratios of the $c_k$, the pieces of the likelihood. (We can assume, without loss of generality, that $\sum_k c_k = 1$.)

We would typically expect the sign of the derivative at $m = 0$ to be positive. Otherwise one of two cases could occur. If the derivative starts out negative and never changes sign, it will be the case that $\hat{m} = 0$ which implies that the Dirichlet model collapses back to the base model. If not, the derivative would cross from negative to positive, implying that a solution to the likelihood equations could result in a minimum, and thus we must exercise more care in finding the MLE. In either case, the sign of the derivative at $m = 0$ depends on the values of the likelihoods with partition sizes $k = 1$ and $k = 2$.

We would also like the limiting sign as $m \to \infty$ to be negative, otherwise it could be the case that $\hat{m} = \infty$ or, depending on the sign at 0, there could be multiple interior extrema. If $c_n$ is close to 0, then the limit can be positive. Note also that the first term of the sign of the derivative equals $\frac{(n-1)(n-2)}{2}$, so for the sign of the derivative to negative, $\frac{c_{n-1}}{c_n}$ should be greater than $\frac{(n-1)(n-2)}{2}$, a number that grows rapidly with the sample size. Thus, the sign of the derivative at $m = \infty$ depends on the likelihoods for partition sizes $k = n - 1$ and $k = n$, and the sample size $n$.

As an illustration of the possible shapes of the likelihood function, we consider a number of simple cases for $n = 6$:

$$
\begin{aligned}
c &= (0, 0, 0, 0, 0, 1), & \frac{\partial}{\partial m}\ell(m) \text{ is } \uparrow \text{ so } \hat{m} = \infty; \\[6pt]
c &= (1, 0, 0, 0, 0, 0), & \frac{\partial}{\partial m}\ell(m) \text{ is } \downarrow \text{ so } \hat{m} = 0; \\[6pt]
c &= (1, 1, 1, 1, 0, 0), & \frac{\partial}{\partial m}\ell(m) \text{ is } - + - \text{ so there is a minimum} \\[2pt]
 & & \text{and maximum;} \\[6pt]
c &= (1, 1, 0, 0, 1, 1), & \frac{\partial}{\partial m}\ell(m) \text{ is } - + \text{ so there is a unique minimum;} \\[6pt]
c &= (1, 0, 0, 0, 0, 1), & \frac{\partial}{\partial m}\ell(m) \text{ is } + - \text{ so there is a unique maximum.}
\end{aligned}
\tag{17}
$$

Thus the likelihood function can have a variety of shapes, which are also illustrated in Figure 1. Moreover, we also note that the likelihood of $m$ can be very flat, meaning that there is insensitivity to the value of the MLE. Referring to (15), the MLE actually corresponds to finding the $m$ that solves this equation. Liu (1996) referred to this as the equating of prior and posterior means (where the expectation is taken using the discrete distribution with weights $m^k c_k$) as the right-hand side of this equation can be interpreted as a prior number of clusters. For example, if we define $\kappa$ by

$$
\kappa = \sum_{i=1}^{n} \frac{m}{m + i - 1},
\tag{18}
$$

then $\kappa$ is the expected number of prior clusters. Even though $\hat{m}$ can be quite variable, there is less variability in $\hat{\kappa} = \sum_{i=1}^{n} \frac{\hat{m}}{\hat{m} + i - 1}$.

3.3. *Posterior mode estimates.* Since the likelihood of the precision parameter depends on the likelihoods from the sub-cluster size $k = 1, 2, n - 1$ and $n$, this reflects an insensitivity to the likelihood and the sample size $n$. For example, when $L_k$ are equal for all $k$, $\lim_{m \to 0} \frac{\partial}{\partial m}\ell(m) > 0$, but $\lim_{m \to \infty} \frac{\partial}{\partial m}\ell(m) > 0$. Thus, in this case, we easily get the MLE as $\hat{m} = \infty$.

Given the potential problems with using and MLE for $m$, we consider a prior on $m$ that results in a unique value of the posterior mode. One of the candidates is a gamma distribution with the shape parameter $a$ and scale parameter $b$. Using the prior $g(m) = \frac{1}{\Gamma(a)b^a} m^{a-1} e^{-m/b}$ we have

$$
L(\theta | \mathbf{y}) = \frac{\Gamma(m)}{\Gamma(m + n)} g(m) \sum_{k=1}^{n} m^k \sum_{A \in \mathcal{A}_k} \prod_{j=1}^{k} \Gamma(n_j) f(\mathbf{y} | \theta, A).
\tag{19}
$$

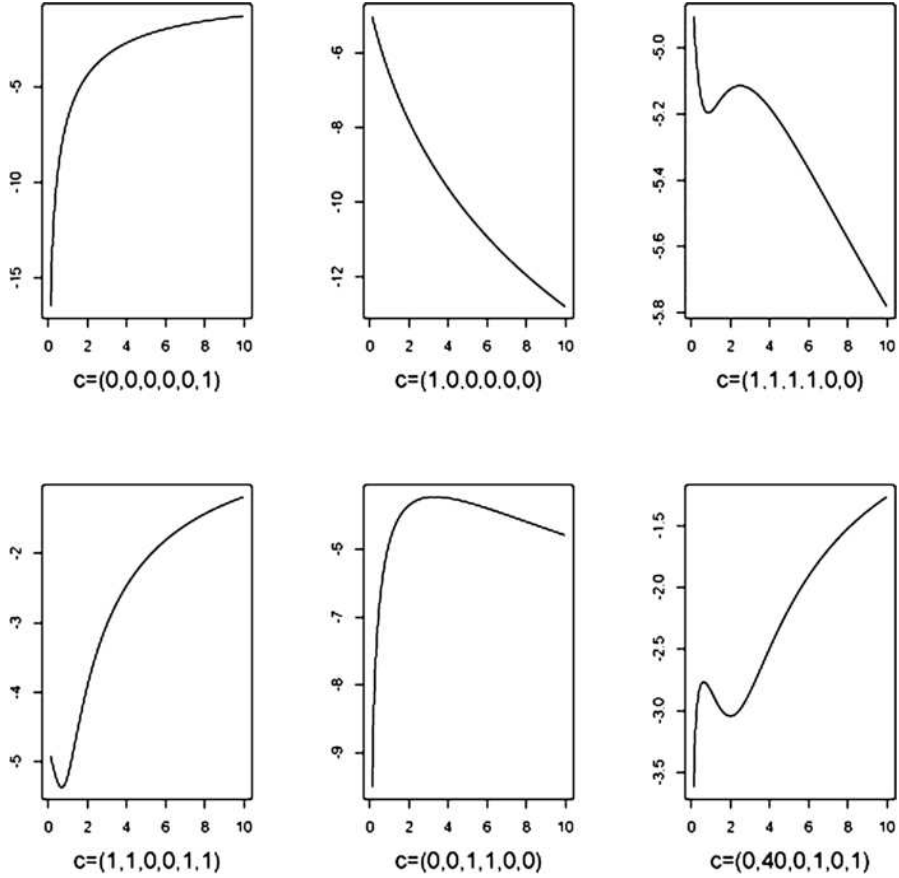

Fig. 1. *Log likelihood functions for a selection of configurations of component likelihoods given in* (17). *The vectors c are the configuration of the marginal likelihoods.*

If we now take logs and differentiate, this amounts to adding the factor $\frac{a-1}{m} - \frac{1}{b}$ to (15). The result of this is that the derivative of this log posterior increases from $m = 0$ and decreases as $m \to \infty$, guaranteeing an interior global maximum. If we integrate (19) we then get the marginal posterior for $m$, which behaves in a similar way.

3.4. *Simulation study of a linear Dirichlet mixed model.* Using the normal linear model of Section 2.2, we conducted a simulation study with a gamma prior, $m \sim \text{gamma}(a, b)$ to study the behavior of the estimates of $m$. We take $n = 6$, $\kappa = 3$, $\mu = 0$, $\tau^2 = 1$ and $\beta = (1, 2, 3)$. The Dirichlet process on the random effect $\boldsymbol{\psi}$ has precision parameter $m$ and base distribution $G_0 = N(0, \tau^2)$. We simulated 100 data sets by generating $X_1$ and $X_2$ from $N(0, 1)$ using the fixed design matrix to generate $Y$.

In this setting, $k = 1, \ldots, 6$, and from the calculation of the Stirling number of the second kind, there are $1, 31, 90, 65, 15, 1$ possible subclusters, respectively. The matrices $A$ associated with these subclusters can be generated, and we sum up the likelihood with all possible subclusters for each $k$ for the profile likelihood.

We estimate $m$ with various settings of the prior mean and variance. With $n = 6$, $\kappa = 3$ and $m = 5$, the solution of $m$ from equation (18) is $m = 1.70$. The numerical summary is given in Table 2 and the histogram of the estimated $k$ is given in Figure 2. For the estimation of $\kappa$, we use the posterior mean of $m$, $\widehat{m}$ and calculate $\widehat{\kappa}$ by using equation (18).

From Table 2, we observe that if the prior mean $ab$ is close to $\widehat{m} = 1.70$, we get a good estimate of $\widehat{\kappa}$ that is close to the fixed $\kappa = 3$. However, if the prior variance is too big, then the estimate of $\kappa$ is less precise. Also, from Figure 2, we observe that the histogram of the estimated $\kappa$ with $ab = 2$ is almost symmetric at $\widehat{\kappa} = 3.10$ with small variance.

## 4. A Gibbs sampler for the model.

We describe a general Gibbs sampling scheme that iteratively generates $A$ matrices and then model parameters assuming that $m$ is fixed at either the MLE or posterior mode. Details on the estimation of $m$ are in Section 6.

Start with the joint likelihood,

$$(20) \qquad L(\theta, A | \mathbf{y}) = \frac{\Gamma(m)}{\Gamma(m+n)} g(m) m^k \prod_{j=1}^{k} \Gamma(n_j) f(\mathbf{y} | \theta, A).$$

With a flat prior on $A$ and $\pi(\theta)$, we get the joint posterior distribution as

$$(21) \qquad \pi(\theta, A | \mathbf{y}) = \frac{m^k f(\mathbf{y} | \theta, A) \pi(\theta)}{\int_{\Theta} \sum_A m^k f(\mathbf{y} | \theta, A) \pi(\theta) \, d\theta}.$$

TABLE 2
*For $n = 6$, $\kappa = 3$, $m = 5$ and various values of the prior parameters, we estimate the precision parameter $m$ and its transformed value $k$. Standard errors are in parentheses*

| Condition | $m$ | $\kappa$ |
|---|---|---|
| $ab = 2$ and $ab^2 = 10$ | 1.85 (0.06) | 3.10 (0.03) |
| $ab = 3$ and $ab^2 = 10$ | 2.56 (0.14) | 3.47 (0.07) |
| $ab = 4$ and $ab^2 = 10$ | 3.07 (0.31) | 3.67 (0.14) |
| $ab = 2$ and $ab^2 = 100$ | 1.98 (0.01) | 3.18 (0.01) |
| $ab = 3$ and $ab^2 = 100$ | 2.95 (0.02) | 3.63 (0.01) |
| $ab = 4$ and $ab^2 = 100$ | 3.90 (0.02) | 3.95 (0.01) |

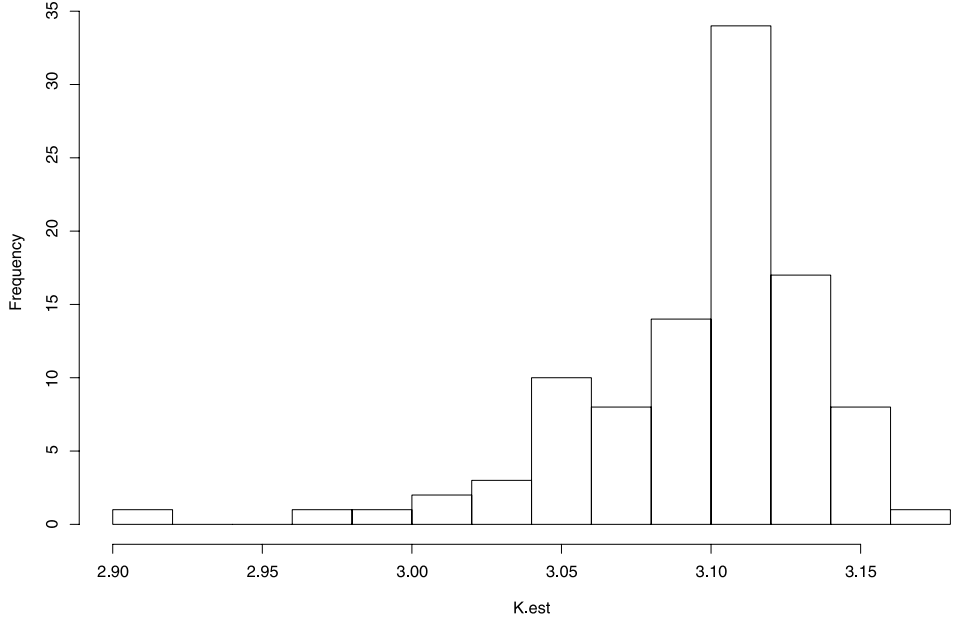

Fig. 2. *A histogram of the estimated k with prior mean* 2 *and variance* 10.

Based on (21), the full conditional posteriors of $\theta$ and $A$ are

$$\pi(\theta|A,\mathbf{y}) = \frac{m^k f(\mathbf{y}|\theta,A)\pi(\theta)}{\int_\Theta m^k f(\mathbf{y}|\theta,A)\pi(\theta)\,d\theta} = \frac{f(\mathbf{y}|\theta,A)\pi(\theta)}{\int_\Theta f(\mathbf{y}|\theta,A)\pi(\theta)\,d\theta},$$

$$\pi(A|\theta,\mathbf{y}) = \frac{m^k f(\mathbf{y}|\theta,A)\pi(\theta)}{\sum_A m^k f(\mathbf{y}|\theta,A)\pi(\theta)} = \frac{m^k f(\mathbf{y}|\theta,A)}{\sum_A m^k f(\mathbf{y}|\theta,A)}.$$

We now outline a Gibbs sampler that will generate from these conditionals by generating $n \times n$ $A$ matrices and recovering the subcluster size through marginalization.

For $t = 1, \ldots, T$, at iteration $t$:

1. Starting from $(\theta^{(t)}, A^{(t)})$,

$$(22) \qquad \theta^{(t+1)} \sim \pi(\theta|A^{(t)}, \mathbf{y}).$$

2. Given $\theta^{(t+1)}$,

$$\mathbf{q}^{(t+1)} = (q_1^{(t+1)}, \ldots, q_n^{(t+1)})$$

$$(23) \qquad \sim \text{Dirichlet}(n_1^{(t)} + r_1, \ldots, n_k^{(t)} + r_k, r_{k+1}, \ldots, r_n),$$

$$A^{(t+1)} \sim m^{k'} f(\mathbf{y}|\theta^{(t+1)}, A) \binom{n}{n_1' \cdots n_n'} \prod_{j=1}^n [q_j^{(t+1)}]^{n_j'},$$

where $n'_1 + \cdots + n'_n = n$ with $k'$ of the $n'_j > 0$. Sampling of the model parameters $\theta$ in (22) is straightforward (Appendix A.1), so we will concentrate on the sampling of $A$ and $\mathbf{q}$.

The matrix $A$ is $n \times n$ with column sums $n_1, \ldots, n_n$, and the columns with zero sums will be removed to obtain an $n \times k'$ matrix, according to Appendix B. Here we keep the $r_j$ as a general choice, but we will see in Section 5.2 and Appendix C.1 that we should choose $r_j = 1$ for all $j$.

The transition kernel of this Markov chain is

$$(24) \qquad k((\theta, A), (\theta', A')) = \pi(\theta'|A, \mathbf{y}) \int_Q P(A'|\mathbf{q}, \theta') f(\mathbf{q}|A) \, d\mathbf{q}$$

with

$$P(A|q, \theta) = \frac{m^k f(\mathbf{y}|\theta, A) \binom{n}{n_1 \cdots n_n} \prod_{j=1}^n q_j^{n_j}}{\sum_A m^k f(\mathbf{y}|\theta, A) \binom{n}{n_1 \cdots n_n} \prod_{j=1}^n q_j^{n_j}}$$

and

$$f(\mathbf{q}|A) = \frac{\Gamma(n + \sum_{j=1}^n r_j)}{\prod_{j=1}^k \Gamma(n_j + r_j) \prod_{j=k+1}^n \Gamma(r_j)} \prod_{j=1}^k q_j^{n_j + r_j - 1} \prod_{j=k+1}^n q_j^{r_j - 1}.$$

Now we take $r_j = 1$, and then we can express the multinomial as

$$\binom{n}{n_1 \cdots n_k} \prod_{j=1}^k q_j^{n_j} = \frac{\Gamma(n+1)}{\prod_{j=1}^n \Gamma(n_j + 1)} \prod_{j=1}^n q_j^{n_j}$$

because the zero valued $n_j$s take care of themselves. With this choice of $r_j$, the transition kernel has $\pi(\theta, A|\mathbf{y})$ as its stationary distribution; details are given in Appendix C.1.

**5. Generating the subclusters.** In this section we discuss two aspects of generating the subclusters. First, we address how to generate according to (23). Then we examine convergence rates and establish that our sampler is an improvement, in terms of operator norm and efficiency, over commonly used algorithms.

5.1. *Generating the matrix $A$.* Generation of the matrix $A$ can be accomplished by using a Gibbs sampler on the rows of $A$. Recall that $a_i, i = 1, \ldots, n$ are the rows of $A$. Define $A_{-i}$ to be the matrix $A$ with the $i$th row removed and $a_i^{(\ell)}$ to be a vector of zeros with a 1 in the $\ell$th position. The matrix $(a_i^{(\ell)}, A_{-i}^{(t)})$ has column sums $n_j^{(\ell)}$ with $k^{(\ell)}$ of $n_j^{(\ell)} > 0$. Then for $i = 1, \ldots, n$,

$$P(a_i = a_i^{(\ell)}|A_{-i}^{(t)})$$

$$\sim \frac{m^{k^{(\ell)}} f(\mathbf{y}|\theta,(a_i^{(\ell)},A_{-i}^{(t)}))\binom{n}{n_1^{(\ell)}\ldots n_n^{(\ell)}}\prod_{j=1}^n [q_j^{(t+1)}]^{n_j^{(\ell)}}}{\sum_{\ell'=1}^n m^{k^{(\ell')}} f(\mathbf{y}|\theta,(a_i^{(\ell')},A_{-i}^{(t)}))\binom{n}{n_1^{(\ell')}\ldots n_n^{(\ell')}}\prod_{j=1}^n [q_j^{(t+1)}]^{n_j^{(\ell')}}},$$

where we update $A_{-i}^{(t)}$ in the usual (Gibbs sampling) way.

Alternatively, we can use a Metropolis–Hastings algorithm with a candidate taken from a multinomial/Dirichlet as described in Appendix B.2. Based on the value of the $q_j$ in (23) we generate a candidate $A$ from the multinomial and then remove the columns with column sum 0. That is, generate an $n \times n$ matrix where each row is a multinomial, and the effective dimension of the matrix, the size of the subclusters, are the non-zero column sums. Deleting the columns with column sum zero is a marginalization of the multinomial distribution. The probability of the candidate follows from Appendix B.2, and the Metropolis–Hastings step is then done.

5.2. *Convergence properties.* From (23), we see that given the subclusters, the sampling of the model parameters from $\pi(\theta|A, \mathbf{y})$ is straightforward. Thus, in investigating convergence we are only concerned with the convergence of the Markov chain on the subclusters. Clearly, if convergence is improved for this part of the chain it will then transfer to the entire chain.

If we ignore the model parameters, then we are concerned with convergence of the chain to the stationary distribution (2), that is,

$$(25) \qquad \pi(A) = \pi(n_1,\ldots,n_k) = \frac{\Gamma(n)}{\Gamma(n+m)} m^k \prod_{j=1}^k \Gamma(n_j),$$

and first we derive the full conditionals in the following way. Start with $(n_1,\ldots,n_k)$ with sum $n-1$, and generate a new row of the $A$ matrix. The matrix $A$ is $n \times k$, and when we generate a new row, either the dimension will remain $n \times k$ or we will increase to $n \times k+1$. If we write $a = \{a_j\}$, then

$$P(a_j=1, n_1,\ldots,n_k) = \frac{\Gamma(n)}{\Gamma(n+m)} m^k \Gamma(n_j+1)\prod_{\substack{j'=1\\j'\neq j}}^k \Gamma(n'_{j'}) \qquad \text{for } j=1,\ldots,k,$$

$$P(a_j=1, n_1,\ldots,n_k) = \frac{\Gamma(n)}{\Gamma(n+m)} m^{k+1}\prod_{j'=1}^k \Gamma(n'_{j'}) \qquad \text{for } j=k+1,$$

$$P(n_1,\ldots,n_k) = \frac{\Gamma(n-1)}{\Gamma(n-1+m)} m^k \prod_{j=1}^k \Gamma(n_j).$$

This results in

$$(26) \qquad P(a_j = 1 | n_1, \ldots, n_k) = \begin{cases} \dfrac{n_j}{n-1+m}, & \text{for } j = 1, \ldots, k, \\ \dfrac{m}{n-1+m}, & \text{for } j = k+1, \end{cases}$$

which are exactly the full conditionals derived by Neal ([2000]), his equation (3.6) ignoring the model parameters, using a limit argument starting from a finite-dimensional Dirichlet. The Gibbs sampler based on (26) is the basis for most of the eight algorithms that he describes; some of which were originally developed by other authors.

The Gibbs sampler resulting from (23), ignoring the model parameters, is

$$(27) \qquad P(A | \mathbf{q}) = \frac{(\Gamma(n)/\Gamma(n+m)) m^k \prod_{j=1}^{k} \Gamma(n_j) \binom{n}{n_1 \cdots n_k} \prod_{j=1}^{k} q_j^{n_j}}{\sum_A (\Gamma(n)/\Gamma(n+m)) m^k \prod_{j=1}^{k} \Gamma(n_j) \binom{n}{n_1 \cdots n_k} \prod_{j=1}^{k} q_j^{n_j}} \quad \text{and}$$

$$f(\mathbf{q} | A) = \frac{\Gamma(n + \sum_{i=1}^{n} r_j)}{\prod_{j=1}^{k} \Gamma(n_j + r_j) \prod_{j=k+1}^{n} \Gamma(r_j)} \prod_{j=1}^{k} q_j^{n_j + r_j - 1} \prod_{j=k+1}^{n} q_j^{r_j - 1},$$

and a similar argument shows that the full conditionals from this chain are

$$(28) \quad P(a_j = 1 | n_1, \ldots, n_k) \propto \begin{cases} \dfrac{n_j}{n_j + 1} \dfrac{q_j}{n-1+m}, & \text{for } j = 1, \ldots, k, \\ \dfrac{m}{n-1+m} q_j, & \text{for } j = k+1, \ldots, n. \end{cases}$$

Notice that for $q_j = n_j + 1$, $j < k$ and $q_j = 1, j > k$ (the normalization does not matter), we see that Neal's Gibbs sampler (26) is the same as (28). We can therefore write the transition kernel of (26) as

$$K_N(A, A') = P(A' | \mathbf{q}^0) g(\mathbf{q}^0 | A),$$

where $g(\mathbf{q}^0 | A)$ is a point mass. In this notation, the kernel of (28) is

$$K(A, A') = \int_Q P(A' | \mathbf{q}) f(\mathbf{q} | \mathbf{q}^0) g(\mathbf{q}^0 | A) \, d\mathbf{q},$$

where $f(\mathbf{q} | \mathbf{q}^0)$ is the same as $f(\mathbf{q} | A)$ in (27). The vector $\mathbf{q}^0$ merely serves to pass the $n_j$.

We are now in the setup of Hobert and Marchev ([2008]) and can use their Theorem 3 to establish the superiority of $K(A, A')$ over $K_N(A, A')$.

THEOREM 1. *For the transition kernels $K_N(A, A')$ and $K(A, A')$, both with stationary distribution given by (25):*

1. *$K(A, A')$ dominates $K_N(A, A')$ in operator norm;*

2. $K(A, A')$ dominates $K_N(A, A')$ in the efficiency ordering of Mira and Geyer (1999) [see also Mira (2001)], which implies that, for any square-integrable function $h$, the asymptotic variance is smaller using $K(A, A')$ than using $K_N(A, A')$.

PROOF. In the terminology of Hobert and Marchev (2008), $K_N(A, A')$ is in the form of a Data Augmentation (DA) algorithm, and $K(A, A')$ is a parameter-expanded version of $K_N(A, A')$. The theorem will be established if we can show that $K(A, A')$ is reversible. This is straightforward as $K(A, A')$ is, itself, a DA algorithm since $K(A, A') = \int_Q P(A'|\mathbf{q}) f(\mathbf{q}|A) \, d\mathbf{q}$. To be specific, if we take $r_j = 1$, then $K(A, A')$ satisfies the detailed balance condition $K(A, A')\pi(A) = K(A', A)\pi(A')$ (Appendix C.2). $\square$

Therefore, in the estimation of any square integrable function $h$, using (28) will result in a smaller variance than obtained by using (26).

5.3. *Assessing the improvement.* The results of Section 5.2 show that our sampler should mix better than "Stickbreaking" as defined by (26). Although we do not know the amount of potential improvement, the results of Roy and Hobert (2007) suggest that there are substantial gains to be had.

To assess the amount of improvement of the Gibbs sampler, the following simulation study was done. For the linear Dirichlet mixed effects model described in Section 2.2 we simulated four data sets. For $n = 100$ we took $A$ matrices corresponding to 1, 5, 25 and 100 groups in the data and ran 20 Markov chains, each for 500 iterations. At each iteration we calculated the variance of the 20 cumulative means which are displayed in Figure 3.

As can be seen, the improvement over the "Stickbreaking" algorithm can be quite substantial; in most cases we see almost a 50% percent improvement. Although we are not claiming that this will hold in all cases, we have a clear indication that substantial reduction in Monte Carlo variance can be attained.

**6. Importance sampling the precision parameter.** To estimate the precision parameter $m$ we want to work with a marginal likelihood function in the form of (7). Based on the development in the previous sections, we start with the marginal posterior,

$$\pi(m|\mathbf{y}) = \frac{\Gamma(m)}{\Gamma(m+n)} g(m) \sum_{k=1}^n m^k f_k(\mathbf{y}), \tag{29}$$

where

$$f_k(\mathbf{y}) = \int_\Theta \sum_{A \in \mathcal{A}_k} \prod_{j=1}^k \Gamma(n_j) f(\mathbf{y}|\theta, A) \pi(\theta) \, d\theta.$$

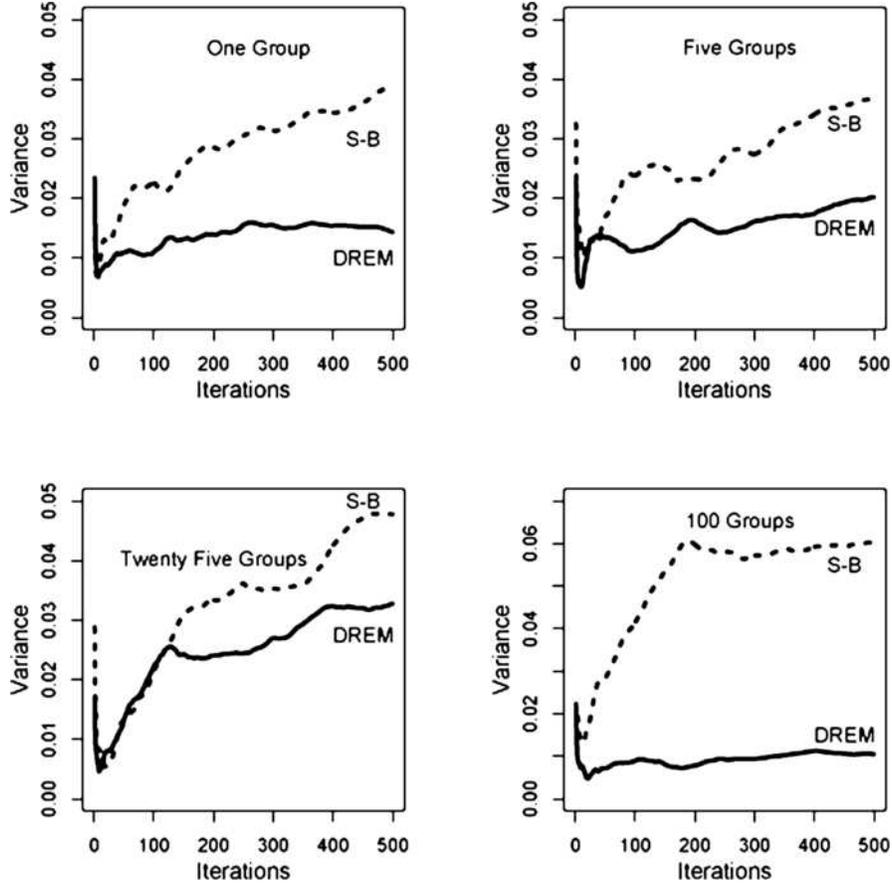

Fig. 3. *For $n = 100$, comparison of variance estimates using the "Stickbreaking" algorithm of (26) (S-B, dashed line) and the algorithm given in (28), the "Dirichlet Random Effects Model" (DREM, solid line). The four plots correspond to four underlying distributions of 1, 5, 25 and 100 groups. Twenty Markov chains were run, and the variance of the 20 estimates was calculated at each of the 500 iterations.*

To take advantage of this functional form for the estimation of $m$, we want to calculate $f_k(\mathbf{y})$ for each $k = 1, \ldots, n$. However, this strategy is difficult to implement for a number of reasons. First, it would necessitate running a full Gibbs sampler (or other MC technique) for all $k = 1, \ldots, n$. Second, the implementation is problematic. For example, consider using an importance sampler based on simulating from the model

$$
\theta \sim f(\mathbf{y}|\theta, A),
$$
$$
(30) \qquad a_i \sim \text{Multinomial}(1, (q_1, \ldots, q_k)), \qquad \text{independent},
$$
$$
\mathbf{q} = (q_1, \ldots, q_k) \sim \text{Dirichlet}(\alpha, \ldots, \alpha),
$$

which leads to the joint posterior distribution

$$(31) \qquad \theta, A \sim f(\mathbf{y}|\theta, A) \frac{\Gamma(k\alpha)}{\prod_{j=1}^{k} \Gamma(\alpha)} \int \prod_{j=1}^{k} q_j^{n_j + \alpha - 1} \, d\mathbf{q}$$

$$= f(\mathbf{y}|\theta, A) \frac{\Gamma(k\alpha)}{\Gamma(n + k\alpha)} \prod_{j=1}^{k} \frac{\Gamma(n_j + \alpha)}{\Gamma(\alpha)},$$

where $n_j = \sum_i a_{ij}$. Unfortunately, there is no guarantee that $n_j > 0$, and samples with $n_j = 0$ will have to be discarded.

However, we can proceed as in (23), and modify (30) to use an $n$-dimensional Dirichlet,

$$\mathbf{q} = (q_1, \dots, q_n) \sim \text{Dirichlet}(\alpha, \dots, \alpha),$$

and then generate $a_i$ independently from this Dirichlet. We then eliminate from the $A$ matrix all columns whose sum is zero. The resulting value of $k$ has the distribution given in (31) because of the marginalization properties of the multinomial and Dirichlet.

The simulation strategy is the following. For $t = 1, \dots, T$ we generate $A^{(t)}$ according to (30) but using the $n$-dimensional Dirichlet, and then marginalize to the number of nonzero $n_j$. We then gather the $A^{(t)}$ according to their values of $k$. Then, for each $k$, if there are $T_k$ matrices $A$ of that size, we estimate $f_k$ by

$$f_k(\mathbf{y}) = \int_\Theta \sum_{A \in \mathcal{A}_k} \prod_{j=1}^{k} \Gamma(n_j) f(\mathbf{y}|\theta, A) \pi(\theta) \, d\theta$$

$$= \sum_{A \in \mathcal{A}_k} \int_\Theta \prod_{j=1}^{k} \frac{\Gamma(n_j) f(\mathbf{y}|\theta, A) \pi(\theta)}{f(\mathbf{y}|\theta, A) \pi(\theta) (\Gamma(k\alpha)/\Gamma(n + k\alpha)) \prod_{j=1}^{k} (\Gamma(n_j + \alpha)/\Gamma(\alpha))}$$

$$(32) \qquad \times f(\mathbf{y}|\theta, A) \pi(\theta) \frac{\Gamma(k\alpha)}{\Gamma(n + k\alpha)} \prod_{j=1}^{k} \frac{\Gamma(n_j + \alpha)}{\Gamma(\alpha)} \, d\theta$$

$$\approx \frac{\Gamma(n + k\alpha)\Gamma(\alpha)^k}{\Gamma(k\alpha)} \frac{1}{T_k} \sum_{t=1}^{T_k} \prod_{j=1}^{k} \frac{\Gamma(n_j^{(t)})}{\Gamma(n_j^{(t)} + \alpha)}$$

$$= \hat{f}_k(\mathbf{y}),$$

where we see very clearly that $m$ only depends on the $n_j$. We now use $\hat{f}_k(\mathbf{y})$ in (29) to obtain the marginal MLE of $m$.

We can further reuse these random variables for all $k' < k$ by randomly choosing two columns and adding them together. This results in an $A$ matrix of one fewer dimension. Details are given in Appendix B.

**7. Application.** In this section we use the Gibbs sampler for a generalized linear mixed model with a Dirichlet process random effect term and probit link to analyze survey data from Scotland. On September 11, 1997, an overwhelming 74.3% of Scottish voters approved of the establishment of the first Scottish national parliament in nearly three hundred years, and on the same ballot the voters gave strong support, 63.5%, to granting this parliament taxation powers. This vote represents a watershed event in the modern history of Scotland which was a free and independent country separate from England until 1707. This vote is part of the Labour government's decentralization program and there is still uncertainty about the future role of Scottish government with the United Kingdom and the European Union. What we are interested in here are those who subsequently voted for the *Conservative* (Tory) party in Scotland and whether such a vote is intended to mitigate Labour's devolution program in Scotland.

The data come from the British General Election Study, Scottish Election Survey, 1997 (ICPSR Study Number 2617). These data contain 880 valid cases, each from an interview with a Scottish national after the election. Our outcome variable of interest is their party choice in the UK general election for Parliament where we collapse all non-Conservative party choices (abstention, Labour, Liberal Democrat, Scottish National, Plaid Cymru, Green, Other, Referendum) to one category which produces 104 Conservative votes. The chosen explanatory variables are intended to explain this choice and include two measures of political efficacy: `POLITICS`, which asks how much interest the respondent has in political events (increasing scale: none at all, not very much, some, quite a lot, a great deal), and `READPAP`, which asks about daily morning reading of the newspapers (yes = 1 or no = 0). It is also important to establish party identity separate from vote choice, `PTYTHNK`, and how strong that party affiliation is for the respondent (categorical by party name), `IDSTRNG` (increasing scale: not very strong, fairly strong, very strong).

We also look at respondents' views on various policy issues. The variable `TAXLESS` asks if "it would be better if everyone paid less tax and had to pay more towards their own healthcare, schools and the like" (measured on a five point increasing Likert scale). `DEATHPEN` asks whether the UK should bring back the death penalty (measured on a five point increasing Likert scale). `LORDS` queries whether the House of Lords should be reformed (asked as *remain as is* coded as zero and *change is needed* coded as one). The question `SCENGBEN` asks how economic benefits are distributed between England and Scotland with the following choices: England benefits more $= -1$, neither/both lose $= 0$, Scotland benefits more $= 1$. The important question `INDPAR` asks which of the following represents the respondent's view on the role of the Scottish government in light of the new parliament: (1) Scotland should become independent, separate from the UK and the

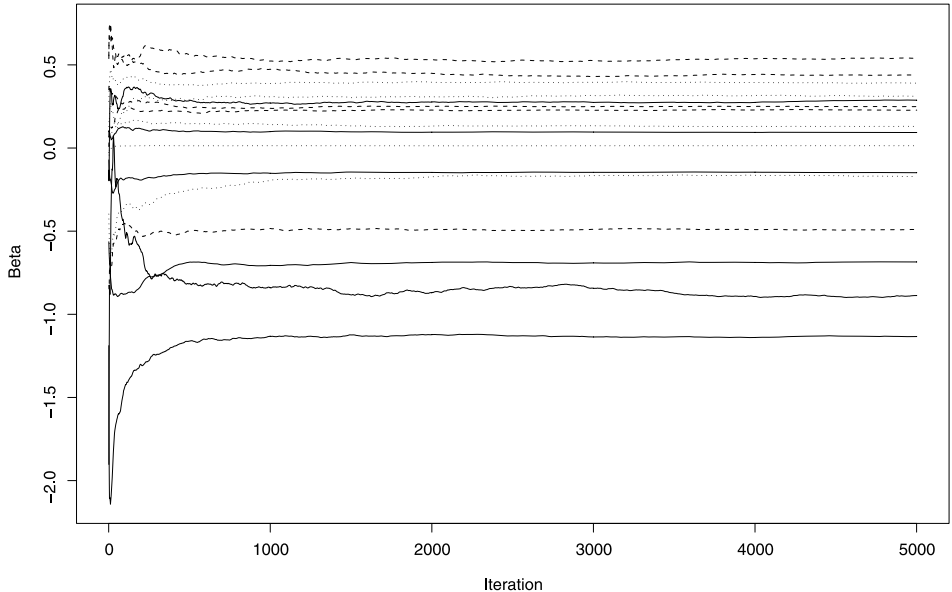

Beta

Iteration

Fig. 4. *Cumulative mean plot, Scotland voting model.*

European Union; (2) Scotland should become independent, separate from the UK but part of the European Union; (3) Scotland should remain part of the UK, with its own elected parliament which has some taxation powers; (4) Scotland should remain part of the UK, with its own elected parliament which has no taxation powers and (5) Scotland should remain part of the UK without an elected parliament. Relatedly, `SCOTPREF1` asks, "should there be a Scottish parliament within the UK?" (yes = 1, no = 0).

Finally, we use three demographic explanatory variables: `RSEX`, the respondent's sex, `RAGE`, the respondent's age, `RSOCCLA2`, the respondents social class (7 category ascending scale), `TENURE1`, whether the respondent rents (0) or owns (1) their household and a categorical variable for church affiliation, measurement of religion is collapsed down to one for the dominant historical religion of Scotland (Church of Scotland/Presbyterian) and zero otherwise and designated `PRESB`.

We set $\sigma^2 = 1$ to establish the scale and provide an intuitive (standard) probit metric. This decision appears to have little influence on the resulting posteriors and allows the $\psi$ specification sufficient latitude to draw nonparametric information from the data. The parameters in the priors on $\mu$ and $\tau^2$ are chosen to make the priors sufficiently diffuse to allow the random effect to do its work. In previous work [Gill and Casella (2009)], we observed little sensitivity to hyperparameter values.

We ran the Gibbs sampler for 5000 iterations disposing of the first 2000. All of the common diagnostics (Geweke, Brooks–Gelman–Rubin, Heidelberger–

TABLE 3
*Posterior model quantiles, voting model*

|  | **0.10** | **0.25** | **0.50** | **0.75** | **0.99** |
|---|---|---|---|---|---|
| CONSTANT | −1.69 | −1.32 | −0.91 | −0.48 | 0.53 |
| POLITICS | 0.12 | 0.17 | 0.23 | 0.29 | 0.44 |
| READPAP | 0.05 | 0.17 | 0.31 | 0.45 | 0.78 |
| PTYTHNK | −0.81 | −0.75 | −0.68 | −0.62 | −0.45 |
| IDSTRNG | 0.14 | 0.19 | 0.25 | 0.31 | 0.46 |
| TAXLESS | 0.02 | 0.07 | 0.13 | 0.19 | 0.32 |
| DEATHPEN | 0.01 | 0.05 | 0.09 | 0.14 | 0.24 |
| LORDS | −0.73 | −0.61 | −0.48 | −0.37 | −0.07 |
| SCENGBEN | 0.22 | 0.30 | 0.39 | 0.47 | 0.67 |
| SCOPREF1 | −1.41 | −1.29 | −1.14 | −1.00 | −0.65 |
| RSEX | 0.20 | 0.30 | 0.43 | 0.56 | 0.84 |
| RAGE | 0.01 | 0.01 | 0.01 | 0.02 | 0.03 |
| RSOCCLA2 | −0.23 | −0.19 | −0.15 | −0.10 | 0.00 |
| TENURE1 | −0.23 | −0.19 | −0.15 | −0.10 | 0.02 |
| PRESB | −0.41 | −0.29 | −0.17 | −0.04 | 0.22 |
| INDPAR | 0.02 | 0.15 | 0.29 | 0.44 | 0.77 |

Welsh, graphics) point toward convergence of the Markov chain to its stationary distribution. Figure 4 is a cumulative mean plot for each of the dimensions for the entire $t = 5000$ period of the chain.

Table 3 provides quantiles for the posterior marginal distributions. We observe that an interest in politics and regular reading of the newspapers increases the probability of voting Conservative as does (not surprisingly) supporting less taxes and the return of the death penalty. We see the same positive effect for men versus women, older versus younger and homeowners versus renters. Those with stronger party attachments are also more likely to vote for the Conservative party. Reforming the House of Lords, Presbyterians and those affiliated with more liberal parties are less likely to vote Conservative.

Two results are surprising. First, those that think that economic policies benefit Scotland more than England are more likely to vote for the Conservative party which is much more aligned with English voters than Scottish voters in general. Perhaps there is a sense that Conservative voting brings attention to Scottish issues from the party. More surprisingly, favoring an independent party is positively associated with voting Conservative through the model. This was our key variable of interest and the relationship is not in the direction expected. The Conservative party is not generally favorable to devolution issues, so these voters are clearly cross-pressured. It is important to keep in mind that the new parliament has taxation powers and thus diminishes the power of local council authorities who are overwhelmingly

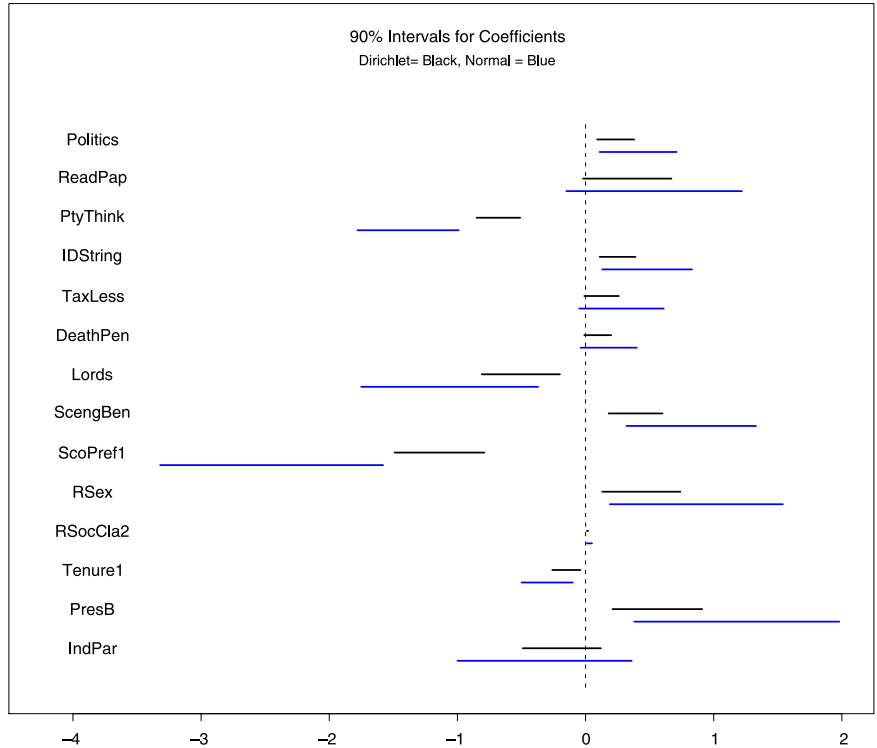

F‌IG. 5. *Comparison of 90% credible intervals for the Dirichlet random effects model (black) to those from a normal random effects model (blue).*

associated with the Labour and Scottish National parties. Thus a Conservative voter may welcome a more centralized taxation program with possibly less influence from these parties, at least at the local level.

In terms of model fit, notice that, aside from the constant, only two marginal posteriors do not have 90% HPD intervals bounded away from zero. It turns out that by every common measure of fit the generalized linear mixed model with a Dirichlet process random effect term outperforms a simple Bayesian probit model with diffuse uniform prior distributions on the parameters. Indeed, when we compare the lengths of credible intervals in Figure 5, we find that the Dirichlet model results in uniformly shorter intervals than those of a normal random effects model. Thus as we anticipated in Section 1, the richer random effects model is able to remove more extraneous variability resulting in tighter credible intervals. We take this as evidence that the new procedure is capturing nonparametric information of interest.

**8. Discussion.** Our interest in models with Dirichlet random effects grew from modeling social science data, where scientists expressed concern over Bayesian models that used informative priors. The Dirichlet process random effects model helps to balance the information from the data and the belief of the researcher while still allowing a normal-type interpretation (in terms of means and variances). As noted previously, we are in agreement with the sentiment of Burr and Doss ([2005](#)), who note that random effects, unlike error terms, can not be checked (there are no residuals). Thus a model with normal random effects is a model of convenience, and moving to a richer model such as the Dirichlet process is a step in relaxing unverifiable assumptions. In particular, the subclustering structure of the Dirichlet process may capture extra variation in the random effects that escape the normal random effects models. The fact that data analysis with the Dirichlet random effect model often differs substantially from the normal random effects model, as noted in Section [7](#), supports this claim.

Representing the subclustering structure through the symmetric binary matrix $A$ is not new. For example, such an equivalence representation was noted by McCullagh and Yang ([2006](#)). Here, the representation has proved useful not only in deriving alternative forms of the model but also in leading to an improved Gibbs sampler. The influence of the random effects, as modeled with the parameters $\psi$ and $\eta$, is only felt through the matrix A, and in some cases these parameters may not have to be generated (see the representation in Section [2.2](#)). This again will lead to a more efficient Gibbs sampler.

The improvement in the Gibbs sampler, as described in Sections [5.2](#) and [5.3](#), appears to come with an increase in computational effort, as we want to start each iteration with an $n \times n$ matrix $A$. However, due to the binary structure of $A$, such a matrix need never be generated. In particular, we can use the correspondence between a multinomial random variable and a discrete random variable to represent the $n \times n$ matrix $A$ as an $n \times 1$ vector $\alpha$. If $X \sim \text{Multinomial}(1, (p_1, \ldots, p_n))$, we create a discrete random variable $X^*$ satisfying $P(X^* = j) = p_j$. We then use $X^*$ to generate the rows of $A$. For example, if $n = 6$ and six samples of $X^*$ gives the vector $(2, 2, 1, 1, 3, 1)$, this represents a matrix $A$ with row 1 having a 1 in column 2, row 2 having a 1 in column 2, etc., with the full matrix being the matrix $A$ in ([5](#)).

We started this project to investigate generalized linear models with Dirichlet random effects but quickly realized that dealing with $m$ is of prime importance and concentrated on linear models to better understand the estimation. As we have seen, ordinary likelihood could be problematic which may be a result of the fact that there is really very little information about $m$ coming from the data. As we saw in Section [3](#), the information in the model about $m$ is only contained in the subclusters, which makes it relatively important to check that the results of the model as somewhat insensitive to the value of the estimated $m$.

Finally, we note that although we have concentrated on linear and probit models, the results will apply directly to a wider class of generalized linear models. There are implementation problems with the Gibbs sampler that arise with models such as the logit, where one needs to use either a slice sampler, a Metropolis–Hastings step or a demarginalization with the Kolmogorov–Smirnov distribution [Andrews and Mallows (1974)]. We have looked at these implementations in Kyung, Gill and Casella (2009). However, these are all variations on the model and, when any general link function such as (9) is used, the Gibbs sampler for $A$ and the estimation of $m$ will remain the same.

## APPENDIX A: GENERATING THE MODEL PARAMETERS

**A.1. A linear model.** For given $A$, the likelihood function is given by

$$L_k(\beta, \sigma^2, \tau^2, \boldsymbol{\eta}|A, \mathbf{y}) = \left(\frac{1}{2\pi\sigma^2}\right)^{n/2} e^{-|\mathbf{y}-X\beta-A\boldsymbol{\eta}|^2/(2\sigma^2)}\left(\frac{1}{2\pi\tau^2}\right)^{k/2} e^{-|\boldsymbol{\eta}|^2/(2\tau^2)}.$$

We add the following normal and inverted gamma (IG) priors:

$$\begin{align}
(33) \qquad \beta|\sigma^2 &\sim N(\mathbf{0}, \sigma^2 I), \\
\tau^2 &\sim \text{IG}(a_1, b_1), \\
\sigma^2 &\sim \text{IG}(a_2, b_2).
\end{align}$$

Then for fixed $m$ and $A$, with $A^* = \frac{1}{\tau^2}I + \frac{1}{\sigma^2}A'A$, a Gibbs sampler of $(\beta, \sigma^2, \tau^2, \boldsymbol{\eta})$ is

$$\boldsymbol{\eta}|\beta, \sigma^2, \tau^2, A, \mathbf{y} \sim N_k\left(\frac{1}{\sigma^2}A^{*-1}A'(\mathbf{y}-X\beta), A^{*-1}\right),$$

$$\beta|\sigma^2, \tau^2, \boldsymbol{\eta}, A, \mathbf{y} \sim N_p((I+X'X)^{-1}X'(\mathbf{y}-A\boldsymbol{\eta}), \sigma^2(I+X'X)^{-1}),$$

$$\tau^2|\beta, \sigma^2, \boldsymbol{\eta}, A, \mathbf{y} \sim \text{IG}\left(\frac{k}{2}+a_1, \frac{1}{2}|\boldsymbol{\eta}|^2+b_1\right),$$

$$\sigma^2|\beta, \tau^2, \boldsymbol{\eta}, A, \mathbf{y} \sim \text{IG}\left(\frac{n+p}{2}+a_2, \frac{1}{2}|\mathbf{y}-X\beta-A\boldsymbol{\eta}|^2+\frac{1}{2}|\beta|^2+b_2\right).$$

If we marginalize out $\boldsymbol{\eta}$, the joint posterior distribution of $(\beta, \sigma^2, \tau^2)$ is

$$\pi_k(\beta, \sigma^2, \tau^2|A, \mathbf{y}) = \int \pi_k(\beta, \sigma^2, \tau^2, \boldsymbol{\eta}|A, \mathbf{y})\, d\boldsymbol{\eta}$$

$$\propto \left(\frac{1}{\sigma^2}\right)^{(n+p)/2+a_2+1}\left(\frac{1}{\tau^2}\right)^{k/2+a_1+1} e^{-(|\beta|^2/2+b_2)/\sigma^2} e^{-b_1/\tau^2}$$

$$\times |A^*|^{1/2} e^{-(\mathbf{y}-X\beta)'[1-A(A^*)^{-1}A'/\sigma^2](\mathbf{y}-X\beta)/(2\sigma^2)}$$

which leads to an alternate Gibbs sampler.

**A.2. A probit model.** Here we need to consider the latent variable $U_i$ such that

$$(34) \qquad U_i = X_i\beta + \psi_i + \eta_i, \qquad \eta_i \sim N(0, \sigma^2),$$

and

$$Y_i = \begin{cases} 1, & \text{if } U_i > 0, \\ 0, & \text{if } U_i \le 0, \end{cases} \qquad i = 1, \ldots, n.$$

It can be shown that $Y_i$ are independent Bernoulli random variables with the probability of success, $p_i = \Phi(\frac{X_i\beta - \psi_i}{\sigma})$ where $\Phi$ is the cdf of the standard Normal.

For given $A$, the likelihood function of model parameters and the latent variable is given by

$$L_k(\beta, \sigma^2, \tau^2, \boldsymbol{\eta}, \mathbf{U}|A, \mathbf{y})$$

$$= \prod_{i=1}^{n}\{I(U_i > 0)I(y_i = 1) + I(U_i \le 0)I(y_i = 0)\}$$

$$\times \left(\frac{1}{2\pi\sigma^2}\right)^{n/2} e^{-|\mathbf{U} - X\beta - A\boldsymbol{\eta}|^2/(2\sigma^2)} \left(\frac{1}{2\pi\tau^2}\right)^{k/2} e^{-|\boldsymbol{\eta}|^2/(2\tau^2)},$$

where $\mathbf{U} = (U_1, \ldots, U_n)$ and leads to the Gibbs sampler

$$\boldsymbol{\eta}|\beta, \sigma^2, \tau^2, \mathbf{U}, A, \mathbf{y} \sim N_k\left(\frac{1}{\sigma^2}A^{*-1}A'(\mathbf{U} - X\beta), A^{*-1}\right),$$

$$\beta|\sigma^2, \tau^2, \boldsymbol{\eta}, A, \mathbf{y} \sim N_p((I + X'X)^{-1}X'(\mathbf{U} - A\boldsymbol{\eta}), \sigma^2(I + X'X)^{-1}),$$

$$\tau^2|\beta, \sigma^2, \boldsymbol{\eta}, A, \mathbf{y} \sim \text{IG}\left(\frac{k}{2} + a_1, \frac{1}{2}|\boldsymbol{\eta}|^2 + b_1\right),$$

$$\sigma^2|\beta, \tau^2, \boldsymbol{\eta}, A, \mathbf{y} \sim \text{IG}\left(\frac{n+p}{2} + a_2, \frac{1}{2}|\mathbf{U} - X\beta - A\boldsymbol{\eta}|^2 + \frac{1}{2}|\beta|^2 + b_2\right),$$

for the model parameters. For the latent variable $\mathbf{U}$, for $i = 1, \ldots, n$,

$$U_i|\beta, \tau^2, \sigma^2, \boldsymbol{\eta}, A, y_i \sim N(X_i\beta + (A\eta)_i, \sigma^2)I(U_i > 0) \qquad \text{if } y_i = 1,$$

$$U_i|\beta, \tau^2, \sigma^2, \boldsymbol{\eta}, A, y_i \sim N(X_i\beta + (A\eta)_i, \sigma^2)I(U_i \le 0) \qquad \text{if } y_i = 0.$$

Here, we can marginalize out $\boldsymbol{\eta}$ such that

$$L_k(\beta, \sigma^2, \tau^2, \mathbf{U}|A, \mathbf{y}) = \int L_k(\beta, \sigma^2, \tau^2, \boldsymbol{\eta}, \mathbf{U}|A, \mathbf{y}) \, d\boldsymbol{\eta}$$

$$= \prod_{i=1}^{n}\{I(U_i > 0)I(y_i = 1) + I(U_i \le 0)I(y_i = 0)\}$$

$$\times \frac{|A^*|^{1/2}}{(2\pi\sigma^2)^{n/2}(\tau^2)^{k/2}}$$

$$\times e^{-(\mathbf{U}-X\beta)'[I-A(A^*)^{-1}A'/\sigma^2](\mathbf{U}-X\beta)/(2\sigma^2)},$$

where $A^* = \frac{1}{\tau^2}I + \frac{1}{\sigma^2}A'A$.

## APPENDIX B: MARGINAL DENSITIES

In this appendix we give the details for the calculation of marginal densities of the Dirichlet, multinomial and their mixture.

**B.1. Dirichlet.** Starting with an $n$-dimensional Dirichlet distribution, the marginal distribution of any $k$ components is also Dirichlet. This corresponds to extracting the rows with non zero column sums in the $A$ matrix in (23). That is, if $(q_1, \ldots, q_n) \sim \text{Dirichlet}(r_1, \ldots, r_n)$, then for $k \leq n$

$$\left( q_1, \ldots, q_{k-1}, \sum_{j=k}^{k+1} q_j, q_{k+2}, \ldots, q_n \right)$$

$$\sim \text{Dirichlet}\left( r_1, \ldots, r_{k-1}, \sum_{j=k}^{k+1} r_j, r_{k+2}, \ldots, r_n \right)$$

as given in Ferguson (1973).

A special case of this result is the combining of two rows which is the marginalization that we use in the calculation of the estimate of $m$ (see Section 6).

If $\mathbf{q} = (q_1, \ldots, q_n) \sim \text{Dirichlet}(r_1, \ldots, r_n)$, then for any $k$ and $k+1 \leq n$

$$\left( q_1, \ldots, q_{k-1}, \sum_{j=k}^{k+1} q_j, q_{k+2}, \ldots, q_n \right)$$

$$\sim \text{Dirichlet}\left( r_1, \ldots, r_{k-1}, \sum_{j=k}^{k+1} r_j, r_{k+2}, \ldots, r_n \right).$$

**B.2. Multinomial.** For $(X_1, \ldots, X_n) \sim \text{Multinomial}(1, (q_1, \ldots, q_n))$, marginalization of the $X_i$s is compatible with the Dirichlet results of the previous section. That is,

$$\left( X_1, \ldots, X_{k-1}, \sum_{j=k}^{n} X_j \right) \sim \text{Multinomial}\left( 1, \left( q_1, \ldots, q_{k-1}, 1 - \sum_{j=1}^{k-1} q_j \right) \right).$$

We also have a similar result for the combining of two elements of the vector, that is, if $(X_1, \ldots, X_n) \sim \text{Multinomial}(1, (q_1, \ldots, q_n))$, then

$$\left( X_1, \ldots, X_{k-1}, \sum_{j=k}^{k+1} X_j, X_{k+2}, \ldots, X_n \right)$$

$$\sim \text{Multinomial} \left( 1, \left( q_1, \ldots, q_{k-1}, 1 - \sum_{j \neq k, k+1} q_j, q_{k+2}, \ldots, q_n \right) \right).$$

**B.3. Multinomial–Dirichlet.** Lastly, we see that these marginalization patterns persist when we combine the multinomial and Dirichlet. Let the matrix $A_{n \times n}$ have each row be an independent multinomial as follows:

$$(a_{i1}, \ldots, a_{in}) \sim \text{Multinomial}(1, (q_1, \ldots, q_n)), \qquad i = 1, \ldots, n,$$

$$(q_1, \ldots, q_n) \sim \text{Dirichlet}(r_1, \ldots, r_n),$$

and then create the matrix $A^*$ by adding together rows $k+1, \ldots, n$, marginalizing $(q_1, \ldots, q_n)$ in the same way. Then

$$P(A^*) = \int P(A|\mathbf{q}) f(\mathbf{q}) \, d\mathbf{q}$$

$$= \frac{\Gamma(\sum_{j=1}^n r_j)}{\prod_{j=1}^n \Gamma(r_j)} \int \prod_{j=1}^{k-1} q_j^{n_j} \left( 1 - \sum_{j=1}^{k-1} q_j \right)^{n_k} \prod_{j=k+1}^n q_j^{r_j-1} \, dq_j,$$

where $n_j = \sum_i a_{ij}, 1, \ldots, k-1$ and $n_k = \sum_i \sum_{j=k}^n a_{ij}$.

If we add together rows $k$ and $k+1$ in the matrix $A$ to obtain $A^*$, a similar result holds:

$$P(A^*) = \frac{\Gamma(\sum_{j=1}^n r_j)}{\prod_{j \neq k, k+1} \Gamma(r_j) \Gamma(\sum_{j=k}^{k+1} r_j)}$$

$$\frac{\prod_{j \neq k, k=1} \Gamma(n_j + r_j) \Gamma(n_k + n_{k+1} + \sum_{j=k}^{k+1} r_j)}{\Gamma(n + \sum_{j=1}^n r_j)}.$$

## APPENDIX C: PROPERTIES OF THE MARKOV CHAIN

**C.1. Stationary distributions of $(\boldsymbol{\theta}, \boldsymbol{A})$.** From the transition kernel of $(\theta, A)$ in (24),

$$\sum_A \int_\Theta K((\theta, A), (\theta', A')) \pi(\theta, A|\mathbf{y}) \, d\theta$$

$$= \sum_A \int_\Theta \pi(\theta'|A, \mathbf{y}) \int_Q P(A'|q, \theta') f(q|A) \, dq \, \pi(\theta, A|\mathbf{y}) \, d\theta$$

$$= \int_Q \left[ \frac{m^k f(\mathbf{y}|\theta', A') \binom{n}{n_1 \cdots n_{k'}} \prod_{j=1}^{k'} q_j^{n'_j}}{\sum_A m^k f(\mathbf{y}|\theta, A) \binom{n}{n_1 \cdots n_k} \prod_{j=1}^{k} q_j^{n_j}} \right]$$

$$\times \sum_A f(q|A) \frac{f(\mathbf{y}|\theta', A) \pi(\theta')}{\int_\Theta f(\mathbf{y}|\theta', A') \pi(\theta') \, d\theta}$$

$$\times \int_\Theta \frac{m^k f(\mathbf{y}|\theta, A) \pi(\theta)}{\sum_A \int_\Theta m^k f(\mathbf{y}|\theta, A) \pi(\theta) \, d\theta} \, d\theta \, dq,$$

and note that the $\Theta$ integral cancels $\int_\Theta f(\mathbf{y}|\theta', A') \pi(\theta') \, d\theta$ inside the sum over $A$. So the integral becomes

$$\int_Q \left[ \frac{m^k f(\mathbf{y}|\theta', A') \binom{n}{n_1 \cdots n_{k'}} \prod_{j=1}^{k'} q_j^{n'_j}}{\sum_A m^k f(\mathbf{y}|\theta, A) \binom{n}{n_1 \cdots n_k} \prod_{j=1}^{k} q_j^{n_j}} \right] \frac{\sum_A m^k f(q|A) f(\mathbf{y}|\theta', A) \pi(\theta')}{\sum_A \int_\Theta f(\mathbf{y}|\theta, A) \pi(\theta) \, d\theta} \, dq.$$

Now take $\beta_j = 1$ for all $j = 1, \ldots, n$ so that

$$f(q|A) = \frac{\Gamma(2n)}{\prod_{j=1}^{k} \Gamma(n_j + 1)} \prod_{j=1}^{k} q_j^{n_j} = \frac{\Gamma(2n)}{n!} \binom{n}{n_1 \cdots n_k} \prod_{j=1}^{k} q_j^{n_j}.$$

This cancels out the denominator sum to leave

$$\frac{\Gamma(2n)}{n!} \frac{\int_Q m^k f(\mathbf{y}|\theta', A') \pi(\theta') \binom{n}{n_1 \cdots n_{k'}} \prod_{j=1}^{k'} q_j^{n'_j} \, dq}{\sum_A \int_\Theta m^k f(\mathbf{y}|\theta, A) \pi(\theta) \, d\theta},$$

and evaluating the integral over $q$ gives

$$\frac{\Gamma(2n)}{n!} \binom{n}{n_1 \cdots n_{k'}} \frac{\prod_{j=1}^{n} \Gamma(n'_j + 1)}{\Gamma(2n)} = 1,$$

where we do the $n$ dimensional integral with $q_j^0$ for $j > k'$. So

$$\sum_A \int_\Theta K((\theta, A), (\theta', A')) \pi(\theta, A|\mathbf{y}) \, d\theta = \frac{m^k f(\mathbf{y}|\theta', A') \pi(\theta')}{\sum_A \int_\Theta m^k f(\mathbf{y}|\theta, A) \pi(\theta) \, d\theta}$$

$$= \pi(\theta', A'|\mathbf{y}).$$

**C.2. Detailed balance.** We have from (27)

$$K(A, A') \pi(A) = \int_Q \left[ \frac{(\Gamma(n)/\Gamma(n+m)) m^{k'} \prod_{j=1}^{k'} \Gamma(n'_j) \binom{n}{n'_1 \cdots n'_k} \prod_{j=1}^{k'} q_j^{n'_j}}{\sum_A (\Gamma(n)/\Gamma(n+m)) m^k \prod_{j=1}^{k} \Gamma(n_j) \binom{n}{n_1 \cdots n_k} \prod_{j=1}^{k} q_j^{n_j}} \right]$$

$$\times \left[ \frac{\Gamma(2n)}{\prod_{j=1}^{k} \Gamma(n_j + 1)} \prod_{j=1}^{k} q_j^{n_j} \right]$$

$$\times \left[ \frac{\Gamma(n)}{\Gamma(n+m)} m^k \prod_{j=1}^{k} \Gamma(n_j) \right] d\mathbf{q}$$

$$= \int_Q \left[ \frac{(\Gamma(n)/\Gamma(n+m)) m^k \prod_{j=1}^{k} \Gamma(n_j) \binom{n}{n_1 \cdots n_k} \prod_{j=1}^{k} q_j^{n_j}}{\sum_A (\Gamma(n)/\Gamma(n+m)) m^k \prod_{j=1}^{k} \Gamma(n_j) \binom{n}{n_1 \cdots n_k} \prod_{j=1}^{k} q_j^{n_j}} \right]$$

$$\times \left[ \frac{\Gamma(2n)}{\prod_{j=1}^{k} \Gamma(n_j' + 1)} \prod_{j=1}^{k'} q_j^{n_j'} \right]$$

$$\times \left[ \frac{\Gamma(n)}{\Gamma(n+m)} m^k \prod_{j=1}^{k'} \Gamma(n_j') \right] d\mathbf{q}$$

$$= K(A', A)\pi(A').$$

M. Kyung
G. Casella
Department of Statistics
University Florida
102 Griffin-Floyd Hall
Gainesville, Florida 32611
USA
E-mail: kyung@stat.ufl.edu
        casella@stat.ufl.edu

J. Gill
Center for Applied Statistics
Washington University
One Brookings Dr., Eliot 315
St. Louis, Missouri
USA
E-mail: jgill@wustl.edu